%% file: CPThurston.tex
\documentclass[dvips,11pt]{article}
\usepackage{amsmath,amssymb,amsfonts,amscd,amsthm,rotating,pstricks,pst-node,bm,xspace}
\input{cmds}

\usepackage[latin1]{inputenc}

\def\bigoplus{\mathop{\oplus}\limits}

\def\prod{\mathop{\Pi}\limits}
\def\sum{\mathop{\Sigma}\limits}

\def\bbf{\boldmath\bfseries}
\def\o{\overline}

\def\Bi{\begin{itemize}}
\def\Ei{\end{itemize}}

\newcommand{\CC}{\mathbb{C}\,}
\newcommand{\RR}{\mathbb{R}\,}
\renewcommand{\rit}{\RR}

\def\vv{{\rm{V}}\,}

\def\Nil{\operatorname{Nil}}
\def\Sol{\operatorname{Sol}}

\def\Repsp{\operatorname{Repsp}}

\def\GL{\operatorname{GL}}

\def\bM{\begin{matrix}}
\def\eM{\end{matrix}}

\def\lr{left (resp. right) }

\def\resp#1{(resp. #1)}
\def\rresp#1{\qquad \mbox{(resp.} \quad #1\ )}

\def\To{\begin{CD} @>>>\end{CD}}

\begin{document}

\setcounter{page}{0}
{\pagestyle{empty}
\include{Thurston00}
\vfill\eject}
\setcounter{page}{1}
\def\thepage{\arabic{page}}
{\parindent=0pt 
\include{Thurston1}
\include{Thurston2}
\include{Thurston3}
\include{ThurstonBib}}
\end{document}

%% file: cmds.tex
\def\simarrow{\mathrel{\raise -0.5mm\hbox{$\sim$}}\hspace{-1.8mm}{\rightarrow} } 

\def\bsimarrow{\leftarrow\hspace{-0.7mm}\mathrel{\raise -0.5mm\hbox{$\backsim$}} }

\def\bt{\begin{tabular}}
\def\te{\end{tabular}}

\def\lettrine#1#2#3{\noindent\hangindent#1\hangafter-#2
\hskip-#1\smash{\hbox to #1{#3\hfill}}\ignorespaces}

\newcommand{\To}[1]{\mathop{\to}\limits_{#1}}

\def\BM{\begin{pmatrix}}
\def\EM{\end{pmatrix}}

\def\ds{\displaystyle}

\def\d=f{\buildrel\hbox{\scriptsize d\'{e}f}\over \Longleftrightarrow}

\def\square{\hfill\hbox{\vrule height .9ex width .8ex depth -.1ex}}
\def\rit{\text{\it I\hskip -2pt  R}}

\def\zit{\text{\it Z\hskip -4pt  Z}}

\def\Bd{{\text B}}

\def\Ed{{\text E}}
\def\Es{{\cal E}}

\def\be{\begin{equation}}
\def\ee{\end{equation}}
\def\beqn{\begin{eqnarray}}
\def\eeqn{\end{eqnarray}}
\def\nobeqn{\begin{eqnarray*}}
\def\noeeqn{\end{eqnarray*}}
\def\ba{\left(\begin{array}}
\def\ea{\end{array} \right) }

\def\epr{\square\vskip 6pt}

\def\u{\underline}
\def\o{\overline}
\def\and{\; \mbox{and} \;}

\def\hfl#1#2{\smash{\mathop{\hbox to 12mm{\rightarrowfill}}
\limits^{\scriptstyle #1}_{\scriptstyle #2}}}

\def\Be{\begin{enumerate}}
\def\Ee{\end{enumerate}}

\def\Bena{\begin{enumerate}
\def\labelenumi{\theenumi)}
\def\theenumi{\arabic{enumi}}
\def\labelenumii{\theenumii)}
\def\theenumii{\alph{enumii}}}

\def\Bean{\begin{enumerate}
\def\labelenumii{\theenumii)}
\def\theenumii{\arabic{enumii}}
\def\labelenumi{\theenumi)}
\def\theenumi{\alph{enumi}}}

\def\Bero{\begin{enumerate}
\def\labelenumii{\theenumii)}
\def\theenumii{\arabic{enumii}}
\def\labelenumi{(\theenumi)}
\def\theenumi{\roman{enumi}}}

\def\BeRo{\begin{enumerate}
\def\labelenumii{\theenumii)}
\def\theenumii{\arabic{enumii}}
\def\labelenumi{(\theenumi)}
\def\theenumi{\Roman{enumi}}}

\def\Bi{\vskip 11pt\begin{itemize}\itemsep=18pt}
\def\bi{\begin{itemize}}
\def\Ei{\end{itemize}\vskip 11pt}
\def\ei{\end{itemize}}

\def\Bd{\begin{description}}
\def\Ed{\end{description}}

\def\R{\right}
\def\L{\left}

%% file: Thurston00.tex
\null\vfill
\begin{center}
{\LARGE The Thurston's program derived from the Langlands global program with singularities}
\vfill
{\sc C. Pierre\/}
\vskip 11pt

Institut de Mathématique pure et appliquée\\
Université de Louvain\\
Chemin du Cyclotron, 2\\
B-1348 Louvain-la-Neuve,  Belgium\\
pierre@math.ucl.ac.be

\vfill

\begin{abstract}
The  seven geometries $H^3$, $S^3$, $H^2\times \rit$, $S^2\times \rit$, ${\rm PSL}_2(\rit)$, $\Nil$ and $\Sol$ of the Thurston's geometrization program are proved to originate naturally from singularization morphisms and versal deformations on euclidean $3$-manifolds generated in the frame of the Langlands global program.

The Poincare conjecture for a $3$-manifold appears as a particular case of this new approach of the Thurston's program.

\end{abstract}
\vfill
\eject
\end{center}

%% file: Thurston1.tex
\section{Introduction}

The classification of $3$-manifolds is now based upon the Thurston's geometrization conjecture \cite{Thu5} which states that every three-dimensional manifold has  a natural decomposition into building blocks characterized by eight specific geometric structures \cite{Pap}.

R. Hamilton \cite{Ham} formalized such approach by introducing the Ricci flow on a Riemannian space and, recently, G. Perelman \cite{Per1}, \cite{Per2}, developed a new method of surgery of the Ricci flow by extending the flow past the singularities in order to solve the technical difficulties of the Hamilton's program \cite{Mor}, \cite{Sch}.
\vskip 11pt

The interesting feature of the Thurston's geometrization conjecture is that it implies the Poincare conjecture in dimension 3 as a special case.
\vskip 11pt

S. Smale \cite{Sma2} proved the Poincare conjecture in the dimensions superior to four by developing a method for breaking manifolds into handles and M. Freedman using Casson-handles succeeded in proving the four-dimensional case \cite{Free}.
\vskip 11pt

But, the three-dimensional case \cite{Sar} still remains difficult, perhaps because it is related to the ``real'' mathematical world \cite{Sma1}.
\vskip 11pt

In this perspective, {\bbf a new approach of the Thurston's program is envisaged in this paper in a ``natural way'', i.e. without envisaging surgeries\/}, in such a way that:
\Bean
\item {\bbf the ``unperturbed'' Thurston's program, corresponding to the Euclidean geometry $E^3$ and to the boundary tori $T^2$, proceeds respectively from the global program of Langlands in dimensions three and two \cite{Pie1}.}

\item {\bbf the ``perturbed'' Thurston's program\/}, associated with the remaining 7 geometries $H^3$, $S^3$, $S^2\times \rit$, $H^2\times \rit$, ${\rm PSL}_2(\rit)$, $\Nil$ and $\Sol$, {\bbf results from the singularities and their versal deformations on $E^3$.}
\Ee
\vskip 11pt

As the structure of the space in dimension three \cite{Hem} corresponds to the real world in which we are living, {\bbf it would be reassuring to know that the unperturbed and perturbed Thurston's programs are linked to the physical reality.\/}  Now, it is the case since it was proved in \cite{Pie3} that the Langlands global program in dimensions two and three with singularities is associated with the generation of the representation spaces of algebraic (bisemi)groups which correspond precisely to the structure of the vacuum space physical fields of elementary particles \cite{Pie4} submitted to strong fluctuations.
\vskip 11pt

In the new context envisaged here, {\bbf the Thurston's program \cite{Thu1}--\cite{Thu4} is directly related to algebraic geometry and number theory\/} by means of the Langlands program {\bbf and to quantum field theory\/} by the geometries of 3-manifolds associated with the field structure of elementary particles.
\vskip 11pt

This new approach of the Thurston's program allows to understand that {\bbf the origin of the geometries of $3$-manifolds \cite{Gor} depends on:\/}
\Bi
\item {\bbf the singularization morphisms\/} which are introduced \cite{Pie2} as the inverse morphisms of desingularizations \cite{Dur}.

\item {\bbf the versal deformations\/} of these singularities.
\item {\bbf the types of geometries characterizing the singularizations and the versal deformations.\/}
\Ei
\vskip 11pt

More concretely, {\bbf the eight geometries of $3$-manifolds will appear as resulting from the following procedure exempt from surgeries:}
\Be
\item {\bbf we start with the Langlands global program in real dimensions two and three\/} by generating the representation spaces $\Repsp(\GL_2(F^T_{\o v}\times F^T_v))$ and $\Repsp(\GL_3(F_{\o v}\times F_v))$ of the algebraic bilinear semigroups $\GL_2(F^T_{\o v}\times F^T_v)$ and $\GL_3(F_{\o v}\times F_v)$ respectively over the products $(F^T_{\o v}\times F^T_v)$ and $(F_{\o v}\times F_v)$ of sets of toroidal and normal real completions.

According to \cite{Pie1}, the representation space
$\Repsp(\GL_2(F^T_{\o v}\times F^T_v))$ is constituted by a tower of conjugacy classes of $\GL_2(F^T_{\o v}\times F^T_v)$ composed of products, right by left, $T^2_R(j)\times T^2_L(j)$ of increasing two-dimensional (semi)tori.  

Similarly, the representation space $\Repsp(\GL_3(F_{\o v}\times F_v))$ is constituted by a tower of conjugacy classes of
$\GL_3(F_{\o v}\times F_v)$ composed of products, right by left, of {\bbf increasing Euclidean semispaces $E^3_R(j)\times E^3_L(j)$, each one characterized by the Euclidean geometry $E^3$.}

This step corresponds to what is called the unperturbed Thurston's program.
\vskip 11pt

\item {\bbf To get the perturbed Thurston's program, we generate singularities by singularization contracting surjective morphisms\/} on the Euclidean semispaces $E^3_R(j)$ and $E^3_L(j)$.

Afterwards, we {\bbf consider the versal deformations of these singularities\/} envisaged as extensions of the singularization morphisms \cite{Pie2}.
\vskip 11pt

\item By analysis of the local geometries resulting from these singularization and versal deformation morphisms, we find that:
\Be
\item {\bbf the singularities of corank $3$, $2$ and $1$ are respectively responsible for the generation of the five local geometries $(H^3$), $(H^2\times \rit,{\rm PSL}_2(\rit))$ and $(\Nil,\Sol)$\/} in the neighborhood of these singularities.

\item {\bbf the versal deformations\/} of these singularities {\bbf in codimensions $3$, $2$ and $1$ generate locally respectively the geometries $S^3$, $S^2\times \rit$ and $\Sol$.}

\item {\bbf the Poincare conjecture for a $3$-manifold appears as a particular case of the versal deformation in codimension $3$.}
\Ee
\Ee

%% file: Thurston2.tex
\section{Origin of the different geometries of 3-manifolds}

Singularization morphisms and versal deformations of these involving non-euclidean geometries will be introduced in this chapter.
\vskip 11pt

But, first, the global Langlands program, leading to the unperturbed Thursdon's program associated with the euclidean geometries, will be recalled.

The developments of this chapter refer to the paper \cite{Pie2}.
\vskip 11pt

\subsection{Representation spaces of algebraic bilinear semigroups}

The Langlands global program is based on the representations of the Weil groups, given by the (representation spaces of) algebraic groups which are in one-to-one correspondence with their cuspidal representations.

So, the representation spaces of algebraic groups have to be taken into account, but, in order to be the most complete and general, bilinear algebraic semigroups will be considered since they cover their linear equivalents \cite{Pie1}.
\vskip 11pt

\Bi
\item {\bbf A bilinear algebraic semigroup $\GL_n(F_R\times F_L)$}:

\Bi
\item is a bilinear semigroup whose bielements are submitted to the cross binary operation $\u\times$ sending products of right and left elements, referring respectively to the lower and upper half spaces, either in diagonal bielements or in cross bielements \cite{Pie1}.

\item can be decomposed according to:
\[ \GL_n(F_R\times F_L)=T^t_n(F_R)\times T_n(F_L)\]
where:
\Bi
\item $F_R$ and $F_L$ are right and left finite algebraic symmetric finite extensions of a global number field $k$ of characteristic zero;
\item $T_n(F_L)$ \resp{$T_n^t(F_R)$} is a \lr semigroup of upper \resp{lower} triangular matrices with entries in the semifield $F_L$ \resp{$F_R$}.
\Ei\Ei

\item {\bbf The representation (bisemi)space of the algebraic bilinear semigroup of matrices $\GL_n(F_R\times F_L)$ is given by the $\GL_n(F_R\times F_L)$-bisemimodule $G^{(n)}(F_R\times F_L)$\/} which is an affine bisemispace $(\vv_R\otimes_{F_R\times F_L}\vv_L)$, in such a way that the \lr affine semispace $\vv_L$ \resp{$\vv_R$} is  localized in the upper \resp{lower} half space, and is symmetric of $\vv_R$ \resp{$\vv_L$}.

\item The left and right equivalence classes of the real completions of $F_L$ and $F_R$ are respectively the infinite real places noted 
$v=\{v _1,\dots,v _j,\dots,v _t\}$ and
$\o v =\{\o v _1,\dots,\o v _j,\dots,\o v _t\}$.

{\bf The infinite places\/} associated with the complex completions of $F_L$ and $F_R$ are the sets 
$w=\{w _1,\dots,w _j,\dots,w _r\}$ and
$\o w =\{\o w _1,\dots,\o w _j,\dots,\o w _r\}$.

The real completions are assumed to cover the corresponding complex completions and the infinite real and complex places are characterized by increasing Galois extension degrees as developed in \cite{Pie2}.

By this way, {\bbf we get a \lr tower\/}
\[
F_v=\{F_{v_1},\dots,F_{v_{j,m_j}},\dots,F_{v_{t,m_t}}\}
\rresp{F_{\o v}=\{F_{\o v_1},\dots,F_{\o v_{j,m_j}},\dots,F_{\o v_{t,m_t}}\}}\]
{\bbf of packets of equivalent real completions covering the \lr tower\/}
\[
F_w=\{F_{w_1},\dots,F_{w_{j}},\dots,F_{w_{t}}\}_{t\equiv r}
\rresp{F_{\o w}=\{F_{\o w_1},\dots,F_{\o w_{j}},\dots,F_{\o w_{t}}\}}\]
{\bf of corresponding complex completions\/}.

\item Let $G^{(n)}(F_{\o v}\times F_v)$ be the bilinear algebraic semigroup with entries in the product, right by left, $F_{\o v}\times F_v$ of towers of packets of equivalent real completions.

$G^{(n)}(F_{\o v}\times F_v)$ is then composed of conjugacy class representatives
$G^{(n)}(F_{\o v_{j,m_j}}\times F_{v_{j,m_j}})$, $1\le j\le t\le \infty $, having multiplicities $1\le m_j \le m^{(j)}$ and being
$\GL_n(F_{\o v_{j,m_j}}\times F_{v_{j,m_j}})$-subbisemimodules $\subset G^{(n)}(F_{\o v}\times F_{v})$.

Similarly, the complex bilinear algebraic semigroup
$G^{(n)}(F_{\o w}\times F_w)$ with entries in the product, right by left, $F_{\o w}\times F_w$ of complex completions, is composed of {\bbf a tower of increasing conjugacy class representatives\/}
$G^{(n)}(F_{\o w_j}\times F_{w_j})$, $1\le j\le r\le \infty $, $t=r$, {\bbf being $\GL_n(F_{\o w_j}\times F_{w_j})$-subbisemimodules covered by the corresponding packets
$G^{(n)}(F_{\o v_j}\times F_{v_j})\}$ of real conjugacy class representatives\/}.
\Ei
\vskip 11pt

\subsection[Proposition: Double tower of increasing $3D$-euclidean subspaces]{\bbf Proposition: Double tower of increasing $3D$-euclidean subspaces}

{\em Let $G_L^{(3)}(F_v)=\Repsp(\GL_3(F_v))\equiv \Repsp(T_3(F_v))$
\resp{$G_R^{(3)}(F_{\o v})=\Repsp(\GL_3(F_{\o v}))\equiv \Repsp(T^t_3(F_{\o v}))$}
denote the representation space of the \lr linear algebraic semigroup of real dimension 3 with entries in $F_v$ \resp{$F_{\o v}$} such that:
\Bi
\item $G_L^{(3)}(F_v)\subset G^{(3)}(F_{\o v}\times F_v)$;
\item $G_R^{(3)}(F_{\o v})\subset G^{(3)}(F_{\o v}\times F_v)$.
\Ei
Then, {\bbf $G_L^{(3)}(F_v)$
\resp{$G_R^{(3)}(F_{\o v})$} is composed of a \lr tower\/}\linebreak
$\{G_L^{(3)}(F_{v_{j,m_j}})\}_{j,m_j}$
\resp{$\{G_R^{(3)}(F_{\o v_{j,m_j}})\}_{j,m_j}$},
$1\le j\le t\le\infty $,
{\bbf of increasing conjugacy class representatives which are three-dimensional euclidean subspaces localized in the upper \resp{lower} half $3D$-space\/}.

Thus we have that:
\[ G_L^{(3)}=E_L^3=\{E_L^3(j,m_j)\}
\rresp{G_R^{(3)}=E_R^3=\{E_R^3(j,m_j)\}}\]
where:
\Bi
\item $E^3_L$ \resp{$E^3_R$} is the $3D$-half upper \resp{lower} euclidean space;
\item $E^3_L(j,m_j)$ \resp{$E^3_R(j,m_j)$} is a $3D$-half upper \resp{lower} euclidean subspace characterized by a rank $r_{E^3_j}\simeq (j\centerdot N)^3$.
\Ei}
\vskip 11pt

\begin{proof}
\Be
\item The conjugacy class representatives
$G_L^{(3)}(F_{v_{j,m_j}})$
\resp{$G_R^{(3)}(F_{\o v_{j,m_j}})$}
are $3D$-half upper \resp{lower} euclidean subspaces $E^3_L(j,m_j)$ \resp{$E^3_R(j,m_j)$} because they are constructed from flat real completions 
$F_{v_{j,m_j}}$
\resp{$F_{\o v_{j,m_j}}$} submitted to the operator ``flat'' morphism \cite{Pie1}
(of a fibre bundle)
\begin{align*}
 T_3(F_{v_{j,m_j}}): \quad F_{v_{j,m_j}}&\To G_L^{(3)}(F_{v_{j,m_j}})\\
\rresp{T^t_3(F_{\o v_{j,m_j}}): \quad F_{\o v_{j,m_j}}&\To G_R^{(3)}(F_{\o v_{j,m_j}})}\end{align*}
sending the completion $F_{v_{j,m_j}}$
\resp{$F_{\o v_{j,m_j}}$} into the \lr 
$\GL_3(F_{v_{j,m_j}})$-subsemimodule $G_L^{(3)}(F_{v_{j,m_j}})$
\resp{$\GL_3(F_{\o v_{j,m_j}})$-subsemimodule $G_R^{(3)}(F_{\o v_{j,m_j}})$}
which is euclidean because {\bbf there is no deviation to euclidicity generated by the injective morphism\/}
$T_3(F_{v_{j,m_j}})$
\resp{$T^t_3(F_{\o v_{j,m_j}})$}
as developed in \cite{Pie2}.

\item $E^3_L(j,m_j)$ and $E_R^3(j,m_j)$ are characterized by a rank $r_{E^3_j}=(j\centerdot N)^3$ because they are built from the completions $F_{v_j}$ and $F_{\o v_j}$ having a Galois extension degree
\[ [F_{v_j}:k]=[F_{\o v_j}:k]=j\centerdot N\]
where $j$ is a global residue degree and $N$ is the degree of an irreducible completion \cite{Pie1}.\qedhere
\Ee
\end{proof}
\vskip 11pt

\subsection[Proposition: Double tower of increasing $2D$-tori]{\bbf Proposition: Double tower of increasing $2D$-tori}

{\em Let $F_v^T$ \resp{$F^T_{\o v}$} be the set of packets of real completions compactified toroidally.

Let $G^{(2)}_L(F^T_v)=\Repsp(T_2(F^T_v))$
\resp{$G^{(2)}_R(F^T_{\o v})=\Repsp(T^t_2(F^T_{\o v}))$}
denote the representation space of the \lr linear algebraic semigroup of real dimension 2 with entries in $F^T_v$ \resp{$F^T_{\o v}$}.

Then, $G_L^{(2)}(F^T_v)$
\resp{$G_R^{(2)}(F^T_{\o v})$}
$\subset G^{(2)}(F^T_{\o v}\times F^T_v)$
is composed of a \lr tower of increasing conjugacy class representatives which are two-dimensional (semi)tori $T^2_L(j)$ \resp{$T^2_R(j)$}, $1\le j\le t\le \infty$, localized in the upper \resp{lower} half $3D$-space.}
\vskip 11pt

\begin{proof}
The toroidal completions $F^T_{v_j}$
\resp{$F^T_{\o v_j}$} of $F^T_v$ \resp{$F^T_{\o v}$} are in fact circles obtained from the corresponding completions $F_{v_j}$ \resp{$F_{\o v_j}$} by connecting their endpoints.

The flat morphism
\begin{align*}
 T_2(F^T_{v_{j}}): \quad F^T_{v_{j}}&\To G_L^{(2)}(F^T_{v_j})\\
\rresp{T_2(F^T_{\o v_{j}}): \quad F^T_{\o v_{j}}&\To G_R^{(2)}(F^T_{\o v_j})}
\end{align*}
sends the toroidal completion $F^T_{v_{j}}$
\resp{$F^T_{\o v_{j}}$} at the $v_j$-th \resp{$\o v_j)$th} place into the conjugacy class representative
$G^{(2)}_L(F^T_{v_{j}})$
\resp{$G^{(2)}_R(F^T_{\o v_{j}})$}
in such a way that
$T^2_L(j)=G_L^{(2)}(F^T_{v_j})$
\resp{$T^2_R(j)=G_R^{(2)}(F^T_{\o v_j})$}.

So,
$T_2(F^T_{v_{j}})$
\resp{$T^t_2(F^T_{\o v_{j}})$}
corresponds to an injective morphism.
\end{proof}
\vskip 11pt

\subsection{Three-dimensional and two-dimensional (semi)sheaves of differentiable functions}

\Bi
\item Let $\phi ^{(3)}_{G_L}(x_{g_L^{(3)}})$
\resp{$\phi ^{(3)}_{G_R}(x_{g_R^{(3)}})$}
denote the set
$\{\phi ^{(3)}_{G_{j,m_{j_L}}}(x_{g_{j_L}^{(3)}})\}_{j,m_j}$
 (resp.\linebreak {$\{\phi ^{(3)}_{G_{j,m_{j_R}}}(x_{g_{j_R}^{(3)}})\}_{j,m_j}$})
of smooth differentiable functions on the set
$\{E^3_L(j,m_j)\}$
(resp.\linebreak  $\{E^3_R(j,m_j)\}$)
of increasing conjugacy class representatives of
$G_L^{(3)}(F_v)$
\resp{$G_R^{(3)}(F_{\o v})$}.  This set
$\phi ^{(3)}_{G_L}(x_{g_L^{(3)}})$
\resp{$\phi ^{(3)}_{G_R}(x_{g_R^{(3)}})$}
of smooth differentiable functions is the {\bbf set of sections of a semisheaf of rings 
$\theta _{G_L^{(3)}}$
\resp{$\theta _{G_R^{(3)}}$}
on the linear algebraic semigroup 
$G_L^{(3)}(F_v)$
\resp{$G_R^{(3)}(F_{\o v})$}} \cite{Pie1}.

\item Similarly,
let $\phi ^{(2)}_{G^T_L}(x_{g_L^{(2)}})$
\resp{$\phi ^{(2)}_{G^T_R}(x_{g_R^{(2)}})$}
denote the set
$\{\phi ^{(2)}_{G_{j_L}}(x_{g_{j_L}^{(2)}})\}_{j}$
(resp.\linebreak {$\{\phi ^{(2)}_{G_{j_R}}(x_{g_{j_R}^{(2)}})\}_{j}$})
of smooth differentiable functions on the set
$\{T^2_L(j)\}_j$
\resp{$\{T^2_R(j)\}_j$}
of increasing two-dimensional (semi)tori of
$G_L^{(2)}(F^T_v)$
\resp{$G_R^{(2)}(F^T_{\o v})$}.  This set
$\phi ^{(2)}_{G^T_L}(x_{g_L^{(2)}})$
\resp{$\phi ^{(2)}_{G^T_R}(x_{g_R^{(2)}})$}
of smooth differentiable functions is the {\bbf set of sections of a semisheaf of rings 
$\theta _{G_L^{(2)}}$
\resp{$\theta _{G_R^{(2)}}$}
on the \lr  algebraic semigroup 
$G_L^{(2)}(F^T_v)$
\resp{$G_R^{(2)}(F^T_{\o v})$}}.
\Ei
\vskip 11pt

\subsection[Proposition: Langlands global correspondence -- $2D$ real case]{\bbf Proposition: Langlands global correspondence -- $2D$ real case}

{\em {\bbf Let $\sigma ^{(2)}_j(W_{F_{\o v_j}}\times W_{F_{v_j}})= G^{(2)}(F_{\o v_j}\times F_{v_j})$ denote the $2$-dimensional representation subspace of the product, right by left, $W_{F_{\o v_j}}\times W_{F_{v_j}}$ of the Weil subgroups restricted to $F_{\o v_j}$ and $F_{v_j}$.}

Let $\phi ^{(2)}_{G_{j_L}}(T^2_L(j))$
\resp{$\phi ^{(2)}_{G_{j_R}}(T^2_R(j))$}, denoting the smooth differentiable \lr function on the (semi)torus $T^2_L(j)$ \resp{$T^2_R(j)$}, be the cuspidal representation 
$\Pi _j(\GL_2(F_{v_j}))$
\resp{$\Pi _j(\GL_2(F_{\o v_j}))$} of the $j$-th conjugacy class representative of the algebraic semigroup $\GL_2(F_v)$ \resp{$\GL_2(F_{\o v})$}.

Then, {\bbf there exists a Langlands global correspondence
\[ T_{L_2}: \quad \sigma ^{(2)}(W^{ab}_{F_R}\times W^{ab}_{F_L})
\To \Pi (\GL_2 ( F_{\o v_\oplus}\times F_{v_\oplus}))\]
between the sum 
$\sigma ^{(2)}(W^{ab}_{F_R}\times W^{ab}_{F_L})$
of the $2$-dimensional conjugacy class representatives of the product, right by left, of the Weil subgroups given by the algebraic bilinear semigroup $G^{(2)}(F_{\o v_\oplus}\times F_{v_\oplus})$ and its cuspidal representation given by 
$\Pi (\GL_2 ( F_{\o v_\oplus}\times F_{v_\oplus}))$, where
$F_{v_\oplus}=\bigoplus\limits_{j,m_j}F_{v_{j,m_j}}$.}}
\vskip 11pt

\begin{proof}
This is rather immediate if we refer to the preprint \cite{Pie1}.

Indeed, the product, right by left, $T^2_R(j)\times T^2_L(j)$ of the (semi)tori $T^2_R(j)$ and $T^2_L(j)$ results from a bijective toroidal compactification of the conjugacy class 
$G^{(2)}(F_{\o v_j}\times F_{v_j})$ of
$G^{(2)}(F_{\o v}\times F_{v})$.

So, we have that: $T^2_R(j)\times T^2_L(j)=G^{(2)}(F^T_{\o v}\times F^T_{v})$.

And, the smooth differentiable bifunction (i.e. the product of a right function by its left equivalent) 
$\phi ^{(2)}_{G_{j_R}}(T^2_R(j))\otimes 
\phi ^{(2)}_{G_{j_L}}(T^2_L(j))$ on
$(T^2_R(j)\times T^2_L(j))$ constitutes the cuspidal representation
$\Pi_j(\GL_2 ( F_{\o v_j}\times F_{v_j}))$
of $\GL_2 ( F_{\o v_j}\times F_{v_j})$
or of $G^{(2)} ( F_{\o v_j}\times F_{v_j})$.

So, the sum $\sigma ^{(2)}(W^{ab}_{F_R}\times W^{ab}_{F_L})$ of the $2$-dimensional conjugacy class representatives of the Weil subgroups given by:
\[
\sigma ^{(2)}(W^{ab}_{F_R}\times W^{ab}_{F_L})
= \bigoplus_j (G^{(2)} ( F_{\o v_j}\times F_{v_j}) )\]
is in one-to-one correspondence with the searched cuspidal representation $\Pi (\GL_2 ( F_{\o v_\oplus}\times F_{v_\oplus}))$ since:
\begin{align*}
\sigma ^{(2)}(W^{ab}_{F_R}\times W^{ab}_{F_L})
& \simeq \bigoplus (\Pi _j (\GL_2 (F_{\o v_j}\times F_{v_j})))\\
& = \Pi (\GL_2 (F_{\o v_\oplus}\times F_{v_\oplus}))\;.
\end{align*}

It then appears that to the set of products, right by left, of the (semi)tori $T^2_R(j)$ and $T^2_L(j)$, $1\le j\le t\le \infty $, referring to the toroidal representation space of the two-dimensional bilinear algebraic semigroup $G^{(2)}(F_{\o v}\times F_v)$, corresponds a cuspidal representation 
$\Pi (\GL_2 (F_{\o v_\oplus}\times F_{v_\oplus}))$, i.e. a Langlands global correspondence.\end{proof}
\vskip 11pt

{\bbf In the three-dimensional case\/}, as we are dealing with the set of products, right by left, of euclidean (semi)spaces $E^3_R(j,m_j)$ and $E^3_L(j,m_j)$, referring to the representation subspaces of the $3D$-bilinear algebraic semigroup $G^{(3)}(F_{\o v}\times F_v)$, 
{\bbf a holomorphic representation of $G^{(3)}(F_{\o v}\times F_v)$ can be built\/} according to \cite{Pie1} but not a cuspidal representation.  Nevertheless, if a toroidal compactification of these subspaces $E^3_R(j,m_j)$ and $E^3_L(j,m_j)$, giving rise to the corresponding $3D$-tori
$T^3_R(j,m_j)$ and $T^3_L(j,m_j)$, is envisaged, then a cuspidal representation of $\GL_3(F_{\o v}\times F_v)$ can be obtained, leading to a $3D$-Langlands global correspondence as it was developed for the two-dimensional case in Proposition 2.5.
\vskip 11pt

\subsection{Introducing singularizations}

The generation of singularities, called {\bf singularizations\/}, will now be recalled; they {\bbf consist of collapses of normal crossings divisors into the singular loci and correspond to contracting surjective morphisms being inverse of those of resolutions of singularities \cite{Dej}, \cite{Hir}, \cite{Zar}.}

These singularizations will be envisaged on the smooth differentiable functions
$\phi ^{(3)}_{G_L}(E^3_L(j,m_j))$
\resp{$\phi ^{(3)}_{G_R}(E^3_R(j,m_j))$}
and $\phi ^{(2)}_{G^T_L}(T^2_L(j))$
\resp{$\phi ^{(2)}_{G^T_R}(T^2_R(j))$} respectively on $3D$-euclidean (semi)spaces and on $2D$-(semi)tori.  To facilitate the notations, they will be written indistinctly $\phi _L$ \resp{$\phi _R$} until the section 2.10.

A normal crossings divisor will be assumed to be a function on one or a set of real irreducible completions of rank $N$ \cite{Pie2}.
\vskip 11pt

\subsection{Proposition: Contracting surjective morphism of singularization}

{\em Let $D_L$ \resp{$D_R$} be a normal crossings divisor of the regular function $\o\phi _L$ \resp{$\o\phi _R$} given by:
\[ \o\phi _L=\phi _L\cup D_L
\rresp{\o\phi _R=\phi _L\cup D_R}.\]
{ \bbf The singularization of $\o\phi _L$ \resp{$\o\phi _R$} into the singular locus $\Sigma _L$ \resp{$\Sigma _R$} results from the contracting surjective morphism:}
\[\o\rho _L: \quad \o\phi _L\To\phi ^*_L
\rresp{\o\rho _R: \quad \o\phi _R\To\phi ^*_R}\]
{\bbf in such a way that:}
\Bean
\item {\bbf $\Sigma _L\subset \phi ^*_L$
\resp{$\Sigma _R\subset \phi ^*_R$} be the union of the homotopic image of 
$D_L\subset \o\phi _L$
\resp{$D_R\subset \o\phi _R$}} and of a possible closed singular sublocus 
$\Sigma ^S_L\subset \Sigma _L$
\resp{$\Sigma ^S_R\subset \Sigma _R$} of $\phi ^*_L$ \resp{$\phi ^*_R$}:
\[ \Sigma _L=\o\rho _L(D_L)\cup\Sigma ^S_L
\rresp{\Sigma _R=\o\rho _R(D_R)\cup\Sigma ^S_R}\;;\]

\item $\o\rho _L$ \resp{$\o\rho _R$} restricted to:
\[ \o\rho ^{\rm is}_L: \quad \phi _L\smallsetminus \o\rho ^{-1}_L(\Sigma ^S_L) \To \phi ^*_L\smallsetminus \Sigma _L
\rresp{\o\rho ^{\rm is}_R: \quad \phi _R\smallsetminus \o\rho ^{-1}_R(\Sigma ^S_R) \To \phi ^*_R\smallsetminus \Sigma _R}\]
be an isomorphism.
\Ee }
\vskip 11pt

\begin{proof}
\Bi
\item The singular sublocus $\Sigma ^S_L\subset \phi ^*_L$
\resp{$\Sigma ^S_R\subset \phi ^*_R$} results from singularizations anterior to that of 
$\o\rho _L(D_L)$
\resp{$\o\rho _R(D_r)$} and then becomes the singular locus of a possible future blowup.

\item Let $\o\rho ^S_L: D_L\to \Sigma _L\smallsetminus \Sigma ^S_L$
\resp{$\o\rho ^S_R: D_R\to \Sigma _R\smallsetminus \Sigma ^S_R$} be the singularization morphism restricted to the singular locus
$\Sigma _L\smallsetminus \Sigma ^S_L$ \resp{$\Sigma _R\smallsetminus \Sigma ^S_R$}.

Then, $\Sigma _L\smallsetminus \Sigma ^S_L$ \resp{$\Sigma _R\smallsetminus \Sigma ^S_R$} is the contracting homotopic image of $D_L$ \resp{$D_R$} in such a way that 
$\o\rho _L\smallsetminus \o\rho ^{\rm is}_L$
\resp{$\o\rho _R\smallsetminus \o\rho ^{\rm is}_R$} be a surjective morphims.

\item {\bf The inverse morphism\/}
\[ \o\rho ^{-1}_L: \quad \phi ^*_L\To\o\phi _L
\rresp{\o\rho ^{-1}_R: \quad \phi ^*_R\To\o\phi _R}\]
{\bf of the singularization $\o\rho _L$ \resp{$\o\rho _R$} corresponds to the blowup of the singular locus} $\Sigma _L$ \resp{$\Sigma _R$} of $\phi ^*_L$ \resp{$\phi ^*_R$} since
\[ \o\rho ^{-1}_L\smallsetminus (\o\rho ^{\rm is}_L)^{-1}: \quad \Sigma _L\To D _L
\rresp{\o\rho ^{-1}_R\smallsetminus (\o\rho ^{\rm is}_R)^{-1}: \quad \Sigma _R\To D _R}\]
is a projective morphism sending the singular locus $\Sigma _L$ \resp{$\Sigma _R$} into the projective normal crossings divisor $D_L$ \resp{$D_R$}.\qedhere
\Ei
\end{proof}
\vskip 11pt

\subsection{Proposition: Sequence of surjective morphisms of singularizations}

{\em
The singularization $\o\rho _L:\o\phi _L\to \phi ^*_L$
\resp{$\o\rho _R:\o\phi _R\to \phi ^*_R$} of the smooth function
$\o\phi _L$
\resp{$\o\phi _R$} is given by the {\bbf following sequence of contracting surjective morphisms\/}:
\begin{align*}
\o\phi _L &\equiv \phi _L^{(0)}
\xrightarrow{\o\rho _L^{(1)}} \phi ^{*(1)}_L
\xrightarrow{\o\rho _L^{(2)}} \phi ^{*(2)}_L
\To\dots
\xrightarrow{\o\rho _L^{(r-1)}} \phi ^{*(r-1)}_L
\xrightarrow{\o\rho _L^{(r)}} \phi ^{*(r)}_L\\
\rresp{\o\phi _R &\equiv \phi _R^{(0)}
\xrightarrow{\o\rho _R^{(1)}} \phi ^{*(1)}_R
\xrightarrow{\o\rho _R^{(2)}} \phi ^{*(2)}_R
\To\dots
\xrightarrow{\o\rho _R^{(r-1)}} \phi ^{*(r-1)}_R
\xrightarrow{\o\rho _R^{(r)}} \phi ^{*(r)}_R}\end{align*}
where $\o\rho ^{(r-1)}_L$ denotes the $(r-1)$-th surjective morphism of singularization of $\o\phi _L$ generating $\phi ^{*(r-1)}_L$, 
in such a way that:
\Bean
\item the singular locus 
$\Sigma _L\subset \phi ^*_L\equiv \phi ^{*(r)}_L$
\resp{$\Sigma _R\subset \phi ^*_R\equiv \phi ^{*(r)}_R$}
is given by:
\begin{align*}
\Sigma _L & \equiv
\Sigma _L^{(r)}=\o\rho _L^{(1)}(D^{(0)}_L) \cup \o\rho _L^{(2)}(D^{(1)}_L) \cup \dots \cup \o\rho ^{(r)}_L(D^{(r-1)}_L)\\
\rresp{\Sigma _R & \equiv
\Sigma _R^{(r)}=\o\rho _R^{(1)}(D^{(0)}_R) \cup \o\rho _R^{(2)}(D^{(1)}_R) \cup \dots \cup \o\rho ^{(r)}_R(D^{(r-1)}_R)}
\end{align*}
where $\Sigma _L^{(1)}=\o\rho _L^{(1)}(D^{(0)}_L)$
\resp{$\Sigma _R^{(1)}=\o\rho _R^{(1)}(D^{(0)}_R)$};

\item $\o\rho _L$ \resp{$\o\rho _R$} restricted to:
\begin{align*}
\o\rho ^{(\rm is)}_L: \quad \o\phi _L\smallsetminus \o\rho _L^{-1}(\Sigma _L^{(r)}) & \To  \phi ^*_L\smallsetminus \Sigma _L^{(r)}\\
\rresp{\o\rho ^{(\rm is)}_R: \quad \o\phi _R\smallsetminus \o\rho _R^{-1}(\Sigma _R^{(r)}) & \To  \phi ^*_R\smallsetminus \Sigma _R^{(r)}}
\end{align*}
is an isomorphism;

\item The orders ``$\ell$'' of the singular subloci $\Sigma ^{(\ell)}_L$ \resp{$\Sigma ^{(\ell)}_R$} form an increasing sequence from left to right, $1\le \ell\le r$.
\Ee
}
\vskip 11pt

\begin{proof}
This proposition is an evident generalization of proposition 2.7 to a set of successive surjective morphisms of singularization giving rise to a singular locus $\Sigma _L^{(r)}$
\resp{$\Sigma _R^{(r)}$} of order ``$r$''.\end{proof}
\vskip 11pt

\subsection{Definition: Corank of the singular locus}

Let $P(x_L,y_L,z_L)$
\resp{$P(x_R,y_R,z_R)$} be the polynomial characterizing the germ of the singular function $\phi ^*_L$ \resp{$\phi ^*_R$} and also the singular locus $\Sigma _L$ \resp{$\Sigma _R$}.

The number of variables of this polynomial is the corank of the germ of $\phi ^*_L$ \resp{$\phi ^*_R$}.

This {\bbf corank is inferior or equal to $3$} according to \cite{A-V-G}.

If the singular locus $\Sigma _L$ \resp{$\Sigma _R$} is given by a singular point of finite codimension, then {\bf the\/} corresponding {\bbf simple germs of differentiable functions are:}
\begin{align*}
A_k: \quad &P(\ x\ )\ =\ x^{k+1}\;, && k\ge 1\;, \\
D_k: \quad &P(x,y)=x^2y+y^{k-1}\;, && k\ge 4\;, \\
E_6: \quad &P(x,y)=x^3+y^4\;,  \\
E_7: \quad &P(x,y)=x^3+xy^3\;,  \\
E_8: \quad &P(x,y)=x^3+y^5\;.\end{align*}
\vskip 11pt

\subsection{The Malgrange division theorem for germs of corank 1}

The Malgrange division theorem for differentiable functions, being the corner stone of the versal deformation, will now be recalled for {\bbf germs of functions having a singularity of corank $1$ and order $k$}.

Let $x'_L=(x_{1_L},x_{2_L},\omega _L)$
\resp{$x'_R=(x_{1_R},x_{2_R},\omega _R)$} be a triple of coordinates in such a way that $x'_L$ \resp{$x'_R$} be localized in the upper \resp{lower} half $3D$-space.

A germ $P(\omega _L)$ \resp{$P(\omega _R)$} has a singularity of corank $1$ and order $k$ in $\omega _L$ \resp{$\omega _R$} if 
$P(0,\omega _L)=\omega ^k_L\ e(\omega _L)$
\resp{$P(0,\omega _R)=\omega ^k_R\ e(\omega _R)$}
where $e(\omega _L)$ \resp{$e(\omega _R)$} is a differentiable unit verifying $e(0)\neq0$.

Let $\theta [\omega _L]$ \resp{$\theta [\omega _R]$} be the algebra of polynomials in $\omega _L$ \resp{$\omega _R$} with coefficients $a(x_L)$ \resp{$a(x_R)$} being ideals of functions defined on a domain $D_{B_L}$ \resp{$D_{B_R}$} included in an open ball centered on $\omega _L$ \resp{$\omega _R$}.

$x_L=(x_{1_L},x_{2_L})$
\resp{$x_R=(x_{1_R},x_{2_R})$} are bituples of coordinates in the upper \resp{lower} half space.

{\bbf The Malgrange division theorem for a germ $P(\omega _L)$ \resp{$P(\omega _R)$} of corank $1$ and order $k$ then corresponds to the versal unfolding of $P(\omega _L)$ \resp{$P(\omega _R)$}} and is given by \cite{Mal}, \cite{Thom}, \cite{Mat}:
\[ f_L=P(\omega _L)\ q_L+R_L
\rresp{f_R=P(\omega _R)\ q_R+R_R}\]
where:
\Bi
\item $f_L$ \resp{$f_R$} is a $3D$-differentiable function (germ);
\item $q_L$ \resp{$q_R$} is a $2D$-differentiable function (germ);
\item $\begin{aligned}[t]
R_L&=\sum\limits^s_{i=1} a_i(x_L)\ \omega ^i_L\in \theta [\omega _L]\\
\rresp{R_R&=\sum\limits^s_{i=1} a_i(x_R)\ \omega ^i_R\in \theta [\omega _R]}\end{aligned}$\\
is a polynomial with degree $s<k$, $s\le3$ in the $3D$-case.
\Ei
The division theorem can be easily stated for germs of corank $2$ and $3$, as developed in \cite{Pie2}.
\vskip 11pt

\subsection{Singularization of semisheaves}

We now come back to the notations of section 2.4 where a \lr semisheaf 
$\theta _{G^{(3)}_L}$
\resp{$\theta _{G^{(3)}_R}$} of $3D$-differentiable functions 
$\phi ^{(3)}_{G_{j,m_{j_L}}}(E^3_L(j,m_j))$
\resp{$\phi ^{(3)}_{G_{j,m_{j_R}}}(E^3_R(j,m_j))$}
on upper \resp{lower} $3D$-euclidean (semi)spaces was introduced.

Similarly, a \lr semisheaf $\theta _{G^{(2)}_L}$
\resp{$\theta _{G^{(2)}_R}$} of $2D$-differentiable functions 
$\phi ^{(2)}_{G_{j_L}}(T^2_L(j))$
\resp{$\phi ^{(2)}_{G_{j_R}}(T^2_R(j))$}
on upper \resp{lower} $2D$-(semi)tori was envisaged.

The two cases will be considered in the following but the developments will only concern here the $3D$-case, the $2D$-case being treated similarly.

{\bbf The singularization of the semisheaf}
$\theta _{G^{(3)}_L}$
\resp{$\theta _{G^{(3)}_R}$} in the sense of proposition 2.8, {\bbf given by the contracting surjective morphism(s):}
\[ \o\rho _{G_L}: \quad \theta _{G^{(3)}_L}\To \theta ^*_{G_L^{(3)}}
\rresp{\o\rho _{G_R}: \quad \theta _{G^{(3)}_R}\To \theta ^*_{G_R^{(3)}}
}\]
{\bbf concerns all the sections} 
$\phi ^{(3)}_{G_{j,m_{j_L}}}(E^3_L(j,m_j)\subset \theta _{G^{(3)}_L}$
\resp{$\phi ^{(3)}_{G_{j,m_{j_R}}}(E^3_R(j,m_j)\subset \theta _{G^{(3)}_R}$} 
which are affected with germs $P_j(\omega _L)$ \resp{$P_j(\omega _R)$} having degenerate singularities of corank inferior or equal to $3$.
\vskip 11pt

\subsection[Proposition: Versal deformation]{Proposition: Versal deformation \cite{G-K}}

{\em
Let $\theta [\omega _L]$
\resp{$\theta [\omega _R]$} be the algebra of polynomials $R_L$ \resp{$R_R$} of the versal unfolding of the germs $P_j(\omega _L)$ \resp{$P_j(\omega _R)$} of the sections of the semisheaf 
$\theta ^*_{G^{(3)}_L}$
\resp{$\theta ^*_{G^{(3)}_R}$}.

Let $\theta _{S_L}=\{\theta (\omega ^1_L),\dots,\theta (\omega ^i_L),\dots,\theta (\omega ^s_L)\}$
\resp{$\theta _{S_R}=\{\theta (\omega ^1_R),\dots,\theta (\omega ^i_R),\dots,\theta (\omega ^s_R)\}$}
denote the family of (semi)sheaves of monomial functions $\omega ^i_L$ \resp{$\omega ^i_R$} of the polynomials 
$R_L\in\theta [\omega _L]$
\resp{$R_R\in\theta [\omega _R]$}.

Then, {\bbf the versal deformation of the singular semisheaf
$\theta ^*_{G^{(3)}_L}$
\resp{$\theta ^*_{G^{(3)}_R}$}
is given by the contracting fiber bundle:}
\[ D_{S_L}: \quad \theta ^*_{G^{(3)}_L}\times \theta _{S_L} \To \theta ^*_{G^{(3)}_L}
\rresp{D_{S_R}: \quad \theta ^*_{G^{(3)}_R}\times \theta _{S_R} \To \theta ^*_{G^{(3)}_R}}\]
{\bbf in such a way that 
$\theta ^{\rm vers}_{G^{(3)}_L}=\theta ^*_{G^{(3)}_L}\times \theta _{S_L}$
\resp{$\theta ^{\rm vers}_{G^{(3)}_R}=\theta ^*_{G^{(3)}_R}\times \theta _{S_R}$}},
being the versal deformation of $\theta ^*_{G^{(3)}_L}$
\resp{$\theta ^*_{G^{(3)}_R}$}, {\bbf is the total space of the fiber bundle $D_{S_L}$ \resp{$D_{S_R})$}} \cite{G-K}, \cite{Mat}.
}
\vskip 11pt

\begin{proof}
The algebra of polynomials $\theta [\omega _L]$ \resp{$\theta [\omega _R]$} is given by
$\theta [\omega _L]=\theta _{S_L}\times\theta _{a_L}$
\resp{$\theta [\omega _R]=\theta _{S_R}\times\theta _{a_R}$} where $\theta _{a_L}$ \resp{$\theta _{a_R}$} is the sheaf of functions $a(x_L)$ \resp{$a(x_R)$} introduced in section 2.10.

{\bbf This algebra of polynomials is the quotient algebra of the versal deformation\/}: it is the quotient of the algebra $\Es_L$ \resp{$\Es_R$} of function germs (generally given by integer power series) by the graded ideal $I_{P_L}$ \resp{$I_{P_R}$} of germs $P(\omega _L)$ \resp{$P(\omega _R)$} \cite{A-V-G}:
\[ \theta [\omega _L]=\Es_L\big/ I_{P_L}
\rresp{\theta [\omega _R]=\Es_R\big/ I_{P_R}}.\]
{\bbf The quotient algebra $\theta [\omega _L]$ \resp{$\theta [\omega _R]$} is thus finitely generated\/}: it is composed of the polynomials $R_L$ \resp{$R_R$} which generate vector (semi)spaces of dimension ``$s$'' which is the codimention of the versal deformation. So, $\theta [\omega _L]$ \resp{$\theta [\omega _R]$} and, thus, $\theta _{S_L}$ \resp{$\theta _{S_R}$} define the versal deformation of the singular semisheaf 
$\theta ^*_{G^{(3)}_L}$
\resp{$\theta ^*_{G^{(3)}_R}$}.
\end{proof}
\vskip 11pt

\subsection{Proposition: Sequence of versal subdeformations}

{\em
The versal unfolding of a germ of differentiable functions is generated by a sequence of contracting morphisms extending the sequence of contracting surjective morphisms of singularization.
}
\vskip 11pt

\begin{proof}
This proposition was proved in \cite{Pie2} for the versal unfolding of the germ 
$P(\omega _L)=\omega _L^{k+1}$
\resp{$P(\omega _R)=\omega _R^{k+1}$}.
Indeed, it was shown that a sequence of $(k-1)$ contracting fiber bundles, whose fibers are divisors projected in the neighbourhood of the singular germ, is responsible for a sequence of $(k-1)$ versal subdeformations generating finitely (i.e. term by term) the polynomial $R_L$ \resp{$R_R$} of the quotient algebra $\theta [\omega _L]$ \resp{$\theta [\omega _R]$}.  The order of the divisors, projected in the neighbourhood of the singular germ, increases in function of the increase of the dimension of the generated vector sub(semi)spaces of the versal unfolding.  By this way, {\bbf the space around the singularity becomes more and more compact in relation with the increase of the (co)dimension of the versal unfolding.}\end{proof}

%% file: Thurston3.tex
\section[Natural generation of the three-dimensional geometries of Thurston]{Natural generation of the three-dimensional geometries of\linebreak Thurston}

Before proving that the non-euclidean $3D$-geometries proceed from singularization morphisms and versal deformations of these, the origin of the hyperbolic and spherical geometries will be introduced.

\subsection{Left and right actions of the Kleinian group}

\Bi
\item The Kleinian group $G_K$ of $\o{\rit}^n=\rit^n\cup \{\infty \}$ is the group of Möbius transformations of $\o{\rit}^n$ acting discontinuously somewhere in $\rit^n$.

{\bbf The action of the Kleinian group $G_K$\/} can be extended to $\o H^{n+1}=H^{n+1}\cup\o{\rit}^n$ where $H^{n+1}=\{(x_1,\dots,x_{n+1})\in\rit^{n+1}:x_{n+1}>0\}$ is the $(n+1)$-dimensional hyperbolic space:
$G_K$ thus {\bbf acts as a group of isometries of $H^{n+1}$ with the hyperbolic metric\/}.

{\bbf The orbit space $M_{G_K}$\/} of the Kleinian group $G_K$ is defined by: $M_{G_K}=(\o H^{n+1}\smallsetminus L(G_K))/G_K$ where{\bbf $L(G_K)$ is the limit set of $G_K$\/} \cite{Mil3}.

This limit set is the closure of the set of fixed points of non-elliptic elements of $G_K$ \cite{Abi}.  It is a nowhere dense set whose area measure is zero.

{\bbf An ordinary set $\Omega (G_K)$\/} of the Kleinian group $G_K$ is given by $\Omega (G_K)=\o{\rit}^n\smallsetminus L(G_K)$.

Recall that a Möbius transformation $g$ of $\o{\rit}^n$ is loxodomic it it is a transformation of the form $g(x)=\lambda \ \alpha (x)$ where $x\in \rit^n$, $\lambda >1$ and $\alpha \in O(n)$ is the orthogonal group of $\rit^n$.  $g$ is hyperbolic if $\alpha ={\rm id.}$, elliptic if $\lambda =1$ and  parabolic if $g$ has the form $g(x)= \alpha (x)+a$, where $a\in\rit^n\smallsetminus\{0\}$.
\vskip 11pt

\item Similarly, {\bbf \lr Möbius transformations $g_L$ \resp{$g_R$}\/}, $g_L\equiv g$, acting discontinuously in the upper \resp{lower} half space $G^{(n)}(F_v)$ \resp{$G^{(n)}(F_{\o v})$} can be introduced as well as the {\bbf \lr action of the Kleinian group $G_K$\/} on the upper \resp{lower} $(n+1)$-dimensional hyperbolic half space 
$H_L^{n+1}\equiv H^{n+1}$
\resp{$H_R^{n+1}$}.  The \lr orbit space
$M_{G_{K_L}}$
\resp{$M_{G_{K_R}}$} associated with the \lr action of the Kleinian group is defined by:
\[
M_{G_{K_L}}=(\o H_L^{n+1}\smallsetminus L(G_{K_L}))\big/G_{K_L}
\rresp{M_{G_{K_R}}=(\o H_R^{n+1}\smallsetminus L(G_{K_R}))\big/ G_{K_R}}\]
where $L(G_{K_L})$ \resp{$L(G_{K_R})$} is the limit set of the Kleinian group $G_K$ acting on the upper \resp{lower} half space.

{\bbf A \lr ordinary set 
$\Omega (G_{K_L})$
\resp{$\Omega (G_{K_R})$}} of $G_{K_L}$ \resp{$G(K_R)$} is given by:
\[ 
\Omega (G_{K_L})= G^{(n)}(F_v)\smallsetminus L(G_{K_L})
\rresp{\Omega (G_{K_R})= G^{(n)}(F_{\o v})\smallsetminus L(G_{K_R})}
\]
where $G^{(n)}(F_v)$ and $G^{(n)}(F_{\o v})$ are introduced in section 2.1.
\Ei
\vskip 11pt

\subsection{Proposition: Hyperbolic geometry in the neighbourhood of the singular locus}

{\em
Let $\Sigma _L$ \resp{$\Sigma _R$} be the singular locus of a germ of corank inferior or equal to $3$ on the singular function $\phi _L^*$ \resp{$\phi _R^*$}.

Let $D_{\Sigma _L}$ \resp{$D_{\Sigma _R}$} denote the neighbourhood of this singular locus.

Thus, we have that:
\Be
\item {\bbf the limit set $L(G_{K_L})$ \resp{$L(G_{K_R})$} of the Kleinian group $G_{K_L}$ \resp{$G_{K_R}$} corresponds to the singular locus $\Sigma _L$ \resp{$\Sigma _R$}.}

\item {\bbf The ordinary set(s) $\Omega (G_{K_L})$ \resp{$\Omega (G_{K_R})$}, characterized by a hyperbolic geometry, correspond to the neighbourhood $D_{\Sigma _L}$ \resp{$D_{\Sigma _R}$} of the singular locus $\Sigma _L$ \resp{$\Sigma _R$}.}
\Ee
}
\vskip 11pt

\begin{proof}
\Be
\item The limit set $L(G_{K_L})$ \resp{$L(G_{K_R})$} is a nowhere dense set, and, furthermore, it has a measure equal to zero: thus, it must correspond to the singular locus $\Sigma _L$ \resp{$\Sigma _R$}.

\item From section 3.1, it then results that the ordinary set $\Omega (G_{K_L})$ \resp{$\Omega (G_{K_R})$} of the Kleinian group $G_{K_L}$ \resp{$G_{K_R})$} is characterized by a hyperbolic geometry.

It thus corresponds to the neighbourhood $D_{\Sigma _L}$ \resp{$D_{\Sigma _R}$} of the singular locus $\Sigma _L$ \resp{$\Sigma _R$}.

On the other hand, it is clear from the sequence of contracting surjective morphisms of singularizations, developed in propositions 2.7 and 2.8, that the neighbourhood of the singular locus must be affected by a hyperbolic geometry.

Finally, it was proved in proposition 2.3.1 of \cite{Pie2} that there is a deviation to Euclidicity in the neighbourhood $D_{\Sigma _L}$ \resp{$D_{\Sigma _R}$} of the singular locus which leads to consider a non-euclidean hyperbolic space of curvature ``$-K$'' on each stratum of $D_{\Sigma _L}$ \resp{$D_{\Sigma _R}$}.\qedhere
\Ee
\end{proof}
\vskip 11pt

\subsection{Proposition: Versal deformation characterized by a spherical geometry}

{\em
Let $P(\omega _L)$ \resp{$P(\omega _R)$} denote a singular germ of corank inferior or equal to $3$ ad codimension $\le 3$ on the $3D$-differentiable function 
$\phi ^{(3)}_{G_{j,m_{j_L}}}(E^3_L(j,m_j))$
\resp{$\phi ^{(3)}_{G_{j,m_{j_R}}}(E^3_R(j,m_j))$}: so,
$\omega _L=(\omega _{1_L})$ 
or $\omega _L=(\omega _{1_L},\omega _{2_L})$
or $\omega _L=(\omega _{1_L},\omega _{2_L},\omega _{3_L})$
\resp{$\omega _R=(\omega _{1_R})$ 
or $\omega _R=(\omega _{1_R},\omega _{2_R})$
or $\omega _R=(\omega _{1_R},\omega _{2_R},\omega _{3_R})$}.

Let $f_L=P(\omega _L)\ q_L+R_L$
\resp{$f_R=P(\omega _R)\ q_R+R_R$}
denote the versal unfolding of the singular germ $P(\omega _L)$ \resp{$P(\omega _R)$} as described in section 2.10.

Then, {\bbf the stratum 
$D_{f_L\mid\phi _L^{(3)}}$
\resp{$D_{f_R\mid\phi _R^{(3)}}$} of the unfolded function $f_L$ \resp{$f_R$} on 
$\phi ^{(3)}_{G_{j,m_{j_L}}}(E^3_L(j,m_j))$
\resp{$\phi ^{(3)}_{G_{j,m_{j_R}}}(E^3_R(j,m_j))$} is characterized by a spherical geometry\/} except perhaps in the neighbourhood of the singular locus $D_{\Sigma _L}$ \resp{$D_{\Sigma _R}$} of the singular germ $P(\omega _L)$ \resp{$P(\omega _R)$} where the geometry is of hyperbolic type.}
\vskip 11pt

\begin{proof}
We are thus concerned with the union of the functions 
$f_L\cup\phi ^{(3)}_{G_{j,m_{j_L}}}(E^3_L(j,m_j))$
\resp{$f_R\cup\phi ^{(3)}_{G_{j,m_{j_R}}}(E^3_R(j,m_j))$} i.e. with the unfolded function $f_L$ \resp{$f_R$} on the $3D$-dimensional substratum function $\phi ^{(3)}_{G_{j,m_{j_L}}}(E^3_L(j,m_j))$
\resp{$\phi ^{(3)}_{G_{j,m_{j_R}}}(E^3_R(j,m_j))$}.

According to section 2.10 and proposition 2.12, this is equivalent to consider in the neighbourhood of the singular germ $P(\omega _L)$ \resp{$P(\omega _R)$}, i.e. on the functions 
$a_i(x_L)\subset \phi ^{(3)}_{G_{j,m_{j_L}}}(E^3_L(j,m_j))$
\resp{$a_i(x_R)\subset \phi ^{(3)}_{G_{j,m_{j_R}}}(E^3_R(j,m_j))$},
$a_i(x_L)\in R_L$
\resp{$a_i(x_R)\in R_R$}, $1\le i\le 3$, the projection of the monomial functions $\omega ^i_L$ \resp{$\omega ^i_R$} of the polynomials $R_L$ \resp{$R_R$} of the quotient algebra of the versal deformation.

Consequently, the stratum 
$D_{f_L\mid\phi _L^{(3)}}$
\resp{$D_{f_R\mid\phi _R^{(3)}}$} of 
 $P(\omega _L)$ \resp{$P(\omega _R)$} on\linebreak 
$\phi ^{(3)}_{G_{j,m_{j_L}}}(E^3_L(j,m_j))$
\resp{$\phi ^{(3)}_{G_{j,m_{j_R}}}(E^3_R(j,m_j))$}
is overcrowded leading to a deviation of Euclidicity characterized by a positive sectional curvature $+K>0$ and by a spherical geometry as proved in \cite{Pie2}.
\end{proof}
\vskip 11pt

\subsection{From the unperturbed Thurston's program to the perturbed one}

We shall now analyze the possible origin of the eight three-dimensional geometries of {\bbf the Thurs\-ton's program which can be stated as follows:}
\begin{quote}
``{\bbf if $M$ is a (closed) oriented prime $3$-manifold, then there is a finite set of disjoint embedded $2$-tori $T^2(j)$ such that each component of the complement in $M$ of $\cup T^2(j)$ admits a geometric structure in the sense of admitting a complete metric, the (necessarily complete) universal metric cover of which is one of the eight three-dimensional model geometries\/} \cite{Gre}.''
\end{quote}

As indicated precedingly, {\bbf the Thurston's program will be split here into an unperturbed and into a perturbed one\/} in such a way that the perturbed part of the Thurston's program results from deformations of the unperturbed Thurston's program due to singularities.
\vskip 11pt

 So, {\bbf the unperturbed Thurston's program will be assumed to originate from the global Langlands program in real dimensions three and two\/} in such a way that:
\Bean
\item the representation space $\Repsp (\GL_3(F_{\o v}\times F_v))$ of the algebraic bilinear semigroup $\GL_3(F_{\o v}\times F_v)$ over the product, right by left, of sets $F_{\o v}$ and $F_v$ of symmetric completions of finite extensions of a global number field $k$ of characteristic zero generates {\bbf a double tower of increasing (compact) three-dimensional euclidean subspaces $E^3_L(j,m_j)$ and $E^3_R(j,m_j)$, $1\le j\le r\le\infty $, localized respectively in the upper and lower half $3D$-spaces \/} as developed in proposition 2.2;

\item the representation space $\Repsp(\GL_2(F^T_{\o v}\times F^T_v))$ of the two-dimensional algebraic bilinear semigroup 
$\GL_2(F^T_{\o v}\times F^T_v)$ over the product, right by left, of sets
$F^T_{\o v}$ and $ F^T_v$ of symmetric toroidal completions generates {\bbf a double tower of increasing two-dimensional tori $T^2_L(j)$ and $T^2_R(j)$} as indicated in proposition 2.3.
\Ee

The reason of considering two-dimensional algebraic bilinear semigroups in the unperturbed\linebreak Thurston's program is that:
\Bean
\item the generated tori $T^2_L(j)$ and $T^2_R(j)$ are two-dimensional compact manifolds.  If they are isotopic to boundary components of $3D$-manifolds, then, these are said to be geometrically atoroidal \cite{Thu1}, \cite{G-L-T}.

\item these tori $T^2_L(j)$ \resp{$T^2_R(j)$} may be glued together pairwise by diffeomorphisms to obtain a closed $3$-manifold or a $3$-manifold with toral boundary \cite{And}.
\Ee
\vskip 11pt

At this stage, we can take up {\bbf the perturbed Thurston's program which can be summarized in the three main propositions.}
\vskip 11pt

\subsection{Proposition: Local geometries round singular loci on $3$-manifolds}

{\em
Assume that the sections 
$\phi ^{(3)}_{G_{j,m_{j_L}}}(E^3_L(j,m_j))$
\resp{$\phi ^{(3)}_{G_{j,m_{j_R}}}(E^3_R(j,m_j))$} of the semisheaf
$\theta _{G_L^{(3)}}$
\resp{$\theta _{G_R^{(3)}}$}
on the euclidean upper \resp{lower} $3$-semispaces
$E^3_L(j,m_j)$ \resp{$E^3_R(j,m_j)$}
are affected by singularization surjective morphisms in such a way that the coranks of their singular germs are inferior of equal to three.

Then, {\bbf the neighbourhood of the singular loci of these singular germs of corank three, two and one are characterized respectively by the local geometries $H^3$, ($H^2\times\rit$ or ${\rm SL}_2(\rit)$) and (${\rm Nil}$ or ${\rm Sol}$)}.
}
\vskip 11pt

\subsection{Proposition: Local geometries of versal deformations on $3$-manifolds}

{\em
Assume that the semisheaf 
$\theta ^*_{G^{(3)}_L}$
\resp{$\theta ^*_{G^{(3)}_R}$}, of which sections on the upper \resp{lower} $3$-semisubspaces are affected by degenerate singularities of corank inferior or equal to three, is submitted to versal deformations of codimensions inferior or equal to three.

Then, {\bbf the neighbourhoods of the unfolded germs in codimensions three, two and one on the sections of 
$\theta ^*_{G^{(3)}_L}$
\resp{$\theta ^*_{G^{(3)}_R}$} are characterized respectively by the local geometries $S^3$, $S^2\times \rit$ and ${\rm Sol}$.
}}
\vskip 11pt

\subsection{Proposition: The Poincare conjecture resulting from the Thurston's program}

{\em \bbf
Assume that the neighbourhood of an unfolded germ in codimension $3$ on a section of $\theta ^*_{G^{(3)}_L}$
or of $\theta ^*_{G^{(3)}_R}$ is a closed simply connected $3$-(semi)manifold.  Then, it is the $3$-``sphere'' $S^3$.  The elimination of the hypothesis of sphericality leads naturally to the Poincare conjecture}
\vskip 11pt

\subsection[Geometric structure on a $3$-manifold]{\bbf Geometric structure on a $3$-manifold}

\Bi
\item The proofs of propositions 3.5, 3.6 and 3.7 are ``clearly'' based on the generation of the three-dimensional local geometries depending on singularities of corank $1$, $2$ and $3$ and on their versal unfoldings in codimensions $1$, $2$ and $3$.

These proofs will be developed in the following, but, before approaching this question, {\bbf the geometric structure on a manifold\/} will be introduced \cite{B-T}, \cite{C-M}, \cite{Kol}.
\vskip 11pt

\item The considered manifolds (or, more exactly, semimanifolds since they are localized in the upper or in the lower half space) are connected (differentiable) manifolds of dimension $3$ generally without boundary \cite{Sha}.

{\bbf A geometric structure on a manifold $M$} is defined by a locally Riemannian metric given by a positive definite quadratic form.

{\bbf The isotropy group\/} of the geometric structure at a point $x\in M$ is defined as the group $G_x$ of linear automorphisms of the tangent space $T_x\ M$ verifying
\[ T_x\ \phi : \quad T_x\ M \To T_x\ M\]
where $T_x\ \phi $ are the differentials of the local isometries sending $x$ into itself \cite{Bon}.

In this respect, let $G$ be the group acting transitively on $M$ in such a way that the stabilizer $G_x$ of $x\in M$ is compact for the compact open topology.

A complete geometric structure on $M$ defines {\bbf a complete $(X,G)$-structure on $M$}, given by an atlas modelling locally $M$ over $X$, where $X$ is the universal covering of $M$ and where $G$ is the isometry group of $X$ \cite{H-R-S}.

{\bbf A geometry in dimension $3$ consists in a pair $(X,G)$ where $X$ is a connected $3$-manifold on which the group $G$ acts transitively.}

$G$ is generally a Lie group and the consideration of a Lie subgroup $H$ leads to the quotient $G/H$ having dimension $3$.  $H$ must be isomorphic to a closed subgroup of $O(3)$.

Thurston introduced eight geometries $(X,G)$ for which there is at least one finite volume complete $(X,G)$-structure \cite{Bon}.
\vskip 11pt

\item For example, the representation space $\Repsp(\GL_3(F_{\o v}\times F_v))$ of the algebraic bilinear semigroup 
$\GL_3(F_{\o v}\times F_v)$ is composed of a left and a right (semi)manifold of which charts are respectively {\bbf three-dimensional euclidean subspaces $E^3_L(j,m_j)$ and $E^3_R(j,m_j)$.}

Thus, this left or right (semi)manifold, having by hypothesis a curvature being equal to zero, has a geometric structure $(X,G)$ in such a way that {\bbf $X$ is isometric to the Euclidean space $E^3$ and $G$ is the isometry group ${\rm Isom}(E^3)$} described by the exact sequence:
\[ O \To \rit^3\To {\rm Isom}(E^3) \To O(3)\To O\;.\]
$G$ is then a discrete group of isometries of $E^3$ and is torsion free.

If $G$ is a finite extension of $\zit$, $G$ is infinite cyclic and $E^3/G$ is the interior of a solid torus or a solid Klein bottle \cite{Sco}.
\vskip 11pt

\item We are now in a position to approach the proofs of propositions 3.5, 3.6 and 3.7.
\Ei
\vskip 11pt

\subsection{Proof of Proposition 3.5: Local geometries round singular loci}

We have to prove that {\bbf the neighbouring chart of the singular locus of a degenerate singular germ of corank three, two and one is characterized respectively by the local geometry $H^3$, $H^2\times \rit$ or ${\rm PSL}_2(\rit)$ and ${\rm Nil}$ or ${\rm Sol}$.}

We refer to the excellent paper of P. Scott \cite{Sco} for the description of these geometries.

\Bean
\item {\bbf Case of corank $3$: the geometry of $H^3$}

It appears from proposition 3.2 that the ordinary sets $\Omega (G_K)$ in the neighbourhood of the singular locus of a degenerate singular germ are characterized by a hyperbolic geometry.  As the envisaged singularity is of corank $3$, the local geometry round the singular locus is the geometry $H^3$ characterized by a negative curvature which is equal to $-1$ if the metric is rescaled \cite{McM}, \cite{Mil1}, \cite{Mil2}.

In fact, {\bbf the neighbourhood of the singular locus is isometric to the hyperbolic $3$-space} \cite{Bra}
\[ H^3=\{(x_1,x_2,x_3)\in\rit^3, x_3>0\}\;.\]
The group of orientations preserving the isometries of $H^3$ is isomorphic to ${\rm PSL}(2,\CC)$: it is the group of Möbius transformations of $\CC\cup\{\infty \}$ given by maps of the form $z\to\dfrac{az+b}{cz+d}$, where $a$, $b$, $c$, $d\in\CC$ and $ad-cb\neq 0$.

The group of complex matrices $\L(\begin{smallmatrix} \ds a & \ds b \\[6pt] \ds c & \ds d\end{smallmatrix}\R)$ acts on $\rit^3_+$, extending its natural action on $\CC\cup\{\infty \}$.
\pagebreak

\item {\bbf Case of corank $2$:}

\Be
\item {\bbf The geometry of $H^2\times \rit$}:

As we are considering {\bbf the local geometry round a singular germ of corank $2$\/}, the ordinary sets $\Omega (G_K)$ round the singular locus {\bbf must be characterized by the hyperbolic geometry $H^2$\/}.  But, as {\bf the charts $M_c$} of the sections of semisheaves 
$\theta ^*_{G^{(3)}_L}$
and $\theta ^*_{G^{(3)}_R}$ are three-dimensional, they {\bbf must be characterized by a local geometry of type $H^2\times \rit$\/} leading to the natural action of the group $G={\rm Isom}(H^2)\times {\rm Isom}(\rit)$.  Thus, the isometry group of $H^2\times \rit$ is isomorphic to ${\rm Isom}(H^2)\times{\rm Isom}(\rit)$, and the local geometry $H^2\times\rit$ is non-isotropic.

If $G$ is the discrete group of isometries of $H^2\times\rit$ having as quotient the chart $M_c$, then the natural foliation of $H^2\times \rit$ by lines gives $M_c$ the structure of a line bundle over some hyperbolic surface in such  a way that $M_c$ cannot be closed \cite{Sco}, \cite{Zhe}.

\item {\bbf The geometry of ${\rm PSL}_2(\rit)$}

But, {\bbf there exists also a twisted version\/} $H^2\overset{\sim}{\times}\rit$ {\bf of} the local geometry {\bbf $H^2\times \rit$ given by $T^1\ H^2$\/} which is the unit tangent bundle of $H^2$, consisting of all tangent vectors of length $1$ of $H^2$.  Topologically, $H^2\overset{\sim}{\times}\rit$ is homeomorphic to $H_2\times\rit$.

The metric of $H^2$ fixes a metric on $T^1\ H^2$ by taking into account that the tangent space $T^1\ H^2$ at $v\in T^1\ H^2$ splits as the direct sum of a line $L_v$ and of a plane $P_v$, where $L_v$ is the tangent line to the fibre $p ^{-1}(p(v))$ where:
\Bi
\item $P_v$ consists of all infinitesimal parallel translations of $v$ along geodesics passing through the point $p(v)\in H^2$.
\item $p:T^1\ H^2\to H^2$ is the natural projection associating its base point to each $v\in T^1\ H^2$.
\Ei

{\bbf There is a natural identification of $T^1\ H^2$ with ${\rm PSL}_2(\rit)$,} the orientation preserving the isometry group of $H^2$ on $T^1\ H^2$.  As the action of this group is transitive and free, the choice of a base point identifies this group with $T^1\ H^2$.  Every orientation preserving the isometry of $H^2$ is a linear fractional map of the form
\[ z\To \dfrac{az+b}{cz+d}\;, \qquad a,b,c,d\in\rit\; , \quad
ad-bc=1\;,\]
which defines a group isomorphism between the orientation preserving the isometry group of $H^2$ and the matrix group ${\rm PSL}_2(\rit)={\rm SL}_2(\rit)/(\pm{\rm Id})$.

If $G$ is a discrete subgroup of isometries of ${\rm PSL}_2(\rit)$ acting on the chart(s) $M_c$, then the foliation of ${\rm PSL}2(\rit)$ by vertical lines gives $M_c$ the structure of a line bundle over a non-closed surface $H^2$.
\Ee
\pagebreak

\item {\bbf Case of corank $1$:}

\Be
\item {\bbf The geometry of ${\rm Nil}$:}

As ${\rm PSL}_2(\rit)$ is a line bundle over $H^2$, there exists a geometry ${\rm Nil}$ which consists in a line bundle over the non-closed Euclidean plane $E^2$: it is a twisted version $E^2\overset{\sim}{\times} E^1$ {\bbf characterized on $\rit^3$ by the Riemannian metric:
\[ ds^2=dx_1^2+dx^2_2+(dx_3-x_1dx_2)^2\;.\]
As there is a negative deviation to Euclidicity in the third dimension ``$x_3$'', there must exist a singular point of corank $1$ in this dimension ``$x_3$''\/}.  This would correspond to a local geometry $H^1\times E^2$ round the singularity since $H^1\subset \rit\simeq E^1$.

The charts $M_c$ on the sections of semisheaves
$\theta ^*_{G^{(3)}_L}$
\resp{$\theta ^*_{G^{(3)}_R}$} are characterized by a Nil geometry of which isometry group $G$ is nilpotent and given by $(3\times 3)$ real upper \resp{lower} unitriangular matrices of the form
$\L(\begin{smallmatrix} \ds1&\ds a&\ds b\\[6pt] \ds 0&\ds1&\ds c\\[6pt] \ds0&\ds0&\ds1\end{smallmatrix}\R)$
\resp{$\L(\begin{smallmatrix} \ds1&\ds0&\ds0\\[6pt]\ds a&\ds 1&\ds 0\\[6pt] \ds b&\ds c&\ds1\end{smallmatrix}\R)$}.  This leads to the exact sequence \cite{Tho}:
\[ O \To \rit\To {\rm Nil} \To \rit^2\To O\]
where $\rit$ consists of the elements of ${\rm Nil}$ with $a=c=0$.

\item {\bbf The geometry of ${\rm Sol}$:}

There still exists a geometry associated with a singularity of corank $1$.  It is the {${\rm Sol}$ geometry characterized by the Riemannian metric
\[ ds^2=e^{+2x_3}\ dx_1^2+e^{-2x_3}\ dx^2_2+dx^2_3\]}
which is such that the discrete group $G$ of transformations of charts of $\theta ^*_{G_L^{(3)}}$ and of $\theta ^*_{G_R^{(3)}}$ acts according to:
\[ ( x_1,x_2,x_3)\To (\epsilon  \ e^{-c}\ x_1+a,\varepsilon '\ e^c\ x_2+b,x_3+c)\]
where $a,b,c\in\rit$.

This group $G$ is defined as a split extension of $\rit^2$ by $\rit$ according to the exact sequence:
\[O\To \rit^2\To {\rm Sol} \To\rit\To O\]
in such a way that $t$ in $\rit$ acts on $\rit^2$ by the map sending $(x_1,x_2)$ to $(e^t\ x_1,e^{-t}\ x_2)$: this corresponds to a linear isomorphism of $\rit^2$ with determinant one and distinct real eigenvalues. Such a linear map is called a hyperbolic isomorphism of $\rit^2$.  And, ${\rm Sol}/G$ is a bundle over a $1$-dimensional orbifold with fibre $S^1\times\rit$ and base $\rit$.  So, {\bbf the local geometry round the degenerate singularity of corank $1$ would be $H^1\times S^1\times\rit$} since:
\Bi
\item the dimension ``$x_2$'' undergoes a negative deviation to Euclidicity due to the factor $e^{-2x_3}$ of $dx_2^2$ in $ds^2$: this explains the hyperbolic geometry $H^1\subset \rit$ in $H^1\times S^1\times \rit$ and the existence of the base $H^1$ of the bundle ${\rm Sol}/G$.

\item the dimension ``$x_1$'' undergoes a positive deviation to Euclidicity due to the factor $e^{+2x_3}$ of $dx_1^2$ in $ds^2$.
\Ei
This is due to the versal unfolding in this dimension, i.e. of codimension one, of the considered degenerate singularity; this explains the spherical geometry $S^1$ in $H^1\times S^1\times \rit$.

Thus, {\bbf any $3$-dimensional chart $M_c$, endowed with a degenerate singularity of corank $1$ and codimension $1$, has a geometric structure modelled on ${\rm Sol}$ and is characterized by a natural foliation associated with the (sub)geometry $S^1\times \rit$.}\epr
\Ee\Ee
\vskip 11pt

\subsection{Proof of Proposition 3.6: Local geometries of versal deformations}

We have to prove that {\bbf the neighbouring chart of an unfolded germ in codimension three, two or one is characterized respectively by a local geometry $S^3$, $S^2\times \rit$ or ${\rm Sol}$.}

\Bean
\item {\bbf Case of codimension $3$: the geometry of $S^3$:}

It appears from proposition 3.3 that the neighbouring chart of an unfolded degenerate singular germ is characterized in the neighbourhood of the singular locus by a spherical geometry.  As the codimension of the considered versal unfolding is equal to three, {\bbf the local geometry round the singular locus must be the geometry of $S^3$} characterized by a positive curvature which is equal to $+1$ if the metric is rescaled \cite{Ber}.  In fact, the envisaged neighbouring chart is isometric to the unit sphere
\[ S^3=\{(x_1,x_2,x_3,x_4)\in\rit^4;\sum^{4}_{i=1}x_i^2=1\}\]
with the Riemannian metric induced by the Euclidean metric of $\rit^4=E^4$.

The isometry group of $S^3$ is ${\rm Isom}(S^3)$ which contains the orthogonal group $O(4)$.

Let $\phi =S^3\to SO(4)$ be the isometry of $S^3$ sending $x\in S^3$ to $q\ x\ q^{-1}$.  Then the image of $\phi $ lies in the subgroup of $SO(4)$ fixing $1$, which can be identified with $SO(3)$.

Remark that the considered neighbouring chart of the unfolded degenerate singular germ in codimension $3$ is not necessarily closed.

\item {\bbf Case of codimension $2$: The geometry of $S^2\times \rit$} \cite{Whi}:

As we are considering the charts or the manifold $M_c$ referring to the versal unfolding in codimension $2$ of a degenerate singular germ,  (their) its local geometry  has   to be of type $S^2$.  But, as the envisaged chart(s) is (are) tridimensional, its (their) local geometry is (are) $S^2\times \rit$ of which isometry group is ${\rm Isom}(S^2)\times{\rm Isom}(\rit)$.

This manifold (or chart(s)) $M_c$ thus has {\bbf the structure of a bundle over the $1$-dimensional base orbifold $\rit$ with fibre $S^2$ resulting locally from a versal deformation in codimension $2$.}

\item {\bbf Case of codimension one: The geometry of ${\rm Sol}$:}

The chart of the versal unfolding in codimension one of a degenerate singular germ must be characterized by the local spherical geometry $S^1$.  As the corank of the versal unfolding in codimension one of a singular germ cannot be generally superior to one and as the envisaged chart of manifold is three-dimensional, the local geometry of a $3$-chart or a $3$-manifold referring to a versal unfolding in codimension one must be $H^1\times S^1\times\rit$, i.e. the ${\rm Sol}$ geometry is considered in section 3.9, c), 2).\epr
\Ee
\vskip 11pt

\subsection{Proof of Proposition 3.7: The Poincare conjecture resulting from the Thurston's program}

This proposition asserts that {\bbf the manifold $M_c$ (or chart) resulting from the versal unfolding in codimension $3$ of a degenerate singular germ is the $3$-sphere $S^3$ if it is closed and simply connected\/}: this corresponds to a natural {\bbf strong version of the Poincare conjecture in dimension $3$}, i.e. when the considered manifold is generated by versal unfolding from a degenerate singular germ.

The Poincare conjecture then appears as originating from the Thurston's geometrization program in the case of a versal deformation in codimension $3$ leading to a manifold having the geometry of $S^3$.  On the other hand, from the standard hypothesis of the Poincare conjecture, we know that the manifold $M_c$ must be closed and simply connected, i.e. composed of closed curves.

Then, it is clear that this manifold $M_c$ is the $3$-sphere $S^3$, which is composed of closed curves (i.e. circles) having the same length: this corresponds to the strong version of the Poincare conjecture.  {\bbf To reduce this strong version of the normal version of the Poincare conjecture, we have to eliminate the hypothesis of sphericality of $M_c$ as resulting from a versal deformation in codimension $3$.}  That is to say, our manifold $M_c$ can be also characterized locally by a hyperbolic and by an euclidean geometry.

Then, $M_c$, being closed, simply connected with all closed curves having the same length, can be deformed continuously by homeomorphism in order to become homeomorphic to $S^3$. This leads to the Poincare conjecture: ``{\bbf Any closed, simply connected $3$-manifold of which closed curves have the same length is homeomorphic to $S^3$,}'' which differs slightly from the classical Poincare conjecture: ``{\bbf Any closed, simply connected $3$-manifold is homeomorphic to $S^3$}.''\epr